\documentclass[11pt, twoside]{article}
\usepackage{amsmath,amsthm}
\usepackage{amsfonts}
\usepackage{amssymb,latexsym}
\usepackage{enumerate}
\newtheorem{theorem}{Theorem}[section]
\newtheorem{remark}[theorem]{Remark}
\newtheorem{proposition}[theorem]{Proposition}
\newtheorem{lemma}[theorem]{Lemma}
\newtheorem{definition}[theorem]{Definition}

\numberwithin{equation}{section}

\def\n{\noindent}
\def\fr{\frac}
\def\Ga{\Gamma}
\def\de{\delta}
\def\ve{\varepsilon}
\def\va{\varphi}
 
\def\veps{\varepsilon}

\def\o{\omega}
\def\O{\Omega}
\def\bb{\bold B}

\def\al{\alpha}
\def\vv{\Vert}

\def\cn{\mathbb C^n}
\def\om{\omega}
\def\pa{\partial}

\def\R{\mathbb{R}}
\def\bC{\mathbb C}

\def\ga{\gamma}
\def\wed{\wedge}
\def\Om{\Omega}

\def\d{\partial}

\begin{document}
	\setlength{\baselineskip}{18truept}
	\pagestyle{myheadings}
	
	%\markboth{  L. M. Hai*, N.Q. Dieu*, N. V. Phu** and Hoang Thieu Anh***}{Hessian type equations for for $m-\omega-$subharmonic functions on bounded domains in $\mathbb C^n$  }
	\title {Hessian type equations for $m-\omega-$subharmonic functions on bounded domains in $\mathbb C^n$}
	\author{
	Hoang Thieu Anh   *,  Le	Mau Hai  **, Nguyen Quang Dieu ** and Nguyen Van Phu  *** \\
		*Faculty of Basic Sciences, University of Transport and Communications, 3 Cau Giay,\\ Dong Da, Hanoi, Vietnam.\\
		 **Department of Mathematics, Hanoi National University of Education,\\ Hanoi, Vietnam.\\
		*** Faculty of Natural Sciences, Electric Power University,\\ Hanoi,Vietnam;
		\\E-mail: anhht@utc.edu.vn, mauhai@hnue.edu.vn \\ngquang.dieu@hnue.edu.vn and  phunv@epu.edu.vn 
	}

	\date{}
	\maketitle
	
	\renewcommand{\thefootnote}{}
	
	
	
	\renewcommand{\thefootnote}{\arabic{footnote}}
	\setcounter{footnote}{0}

\begin{abstract}
	\n
In this paper, we study Hessian type equations for $m-\omega$ subharmonic functions.
Using the  recent results in \cite{KN23a}, \cite{KN23b}, we are able to show
the existence of bounded solutions for such equations on bounded domains in $\mathbb{C}^n$. 
	\end{abstract}
	
\section{Introduction}
Let \(\Omega \subset \mathbb{C}^n\) be a bounded domain. In their seminal contributions from the early 1980s, E. Bedford and B. A. Taylor \cite{BT1} introduced and investigated the complex Monge-Amp\`ere operator $(dd^c .)^n$ for the class of plurisubharmonic functions. They proved that this operator is well-defined on the class of locally bounded plurisubharmonic (psh) functions, and that it yields non-negative Radon measures. This foundational result paved the way for posing the Dirichlet problem for the complex Monge-Ampère equation with a positive Radon measure $\mu$ on $\Omega$. Namely, let $\O\subset \cn$ be a bounded domain and $\mu$ a positive Radon measure on $\O$. Assume that $\varphi\in C^0(\pa\O)$. The Dirichlet problem is to find 
\[
\begin{cases}
u \in \mathrm{PSH}(\Omega) \cap L^{\infty}_{loc}(\overline{\Omega}), \\
(dd^c u)^n =  \mu, \\
\lim\limits_{\Omega \ni z \to x} u(z) = \varphi(x), \quad \forall x \in \partial \Omega.
\end{cases}
\tag{1.1} \label{bt1.0}
\]

\n  where $PSH(\O)$ denotes the set of plurisubharmonic functions on $\O$. It has been shown by Bedford and Taylor in \cite{BT76} that if $\O\subset\cn$ is a strictly pseudoconvex domain and $\mu = fdV_{2n}, f\in C(\overline{\O})$, where $dV_{2n}$ denotes the Lebesgue measure of $\cn$ then (1.1)  is solvable and the solution $u$ belongs to $C(\overline{\O})$. It is known that continuous solutions also exist for $\mu=fdV_{2n}$ where $f\in L^2(\O,dV_{2n})$ (see \cite{CePe92}).  Next, by using suitable techniques of pluripotential theory, Ko{\l}odziej
in \cite{Ko96, Ko98} has proved that (1.1) admits continuous solutions if $\mu = fdV_{2n}$, where $f\in L^p(\O ,dV_{2n}), p>1$. After that in \cite{K95},  Ko{\l}odziej has shown that if there exists a subsolution for the Dirichlet problem (1.1), then the problem is solvable. Now we deal with the Dirichlet problem for Monge-Amp\`ere type equation, an extension of Monge-Amp\`ere equation (1.1). Let $\O\subset\cn$ be a bounded domain, $\mu$ a positive Radon measure on $\O$. Assume that $F:\mathbb{\R}\times\O\longrightarrow [0,+\infty)$ and $\varphi\in C(\pa\O)$ are given. The Dirichlet problem for Monge-Amp\`ere type equation is to find 
\[
\begin{cases}
u \in \mathrm{PSH}(\Omega) \cap L^{\infty}(\overline{\Omega}), \\
(dd^c u)^n = F(u,z) d\mu, \\
\lim\limits_{\Omega \ni z \to x} u(z) = \varphi(x), \quad \forall x \in \partial \Omega.
\end{cases}
\tag{1.2} \label{bt1.2}
\]

\n This problem was first published by Bedford and Taylor in \cite{BT79} when  $d\mu = dV_{2n}$ and $\Omega$ is a bounded strongly pseudoconvex domain. They demonstrated the existence and uniqueness of solutions under the assumption that $F(t,z)\in C^0(\mathbb{R} \times \overline{\Omega})$ and $F^{1/n}$ is convex and non-decreasing in the first variable. In 1984, Cegrell in \cite{Ce84} extended these results by proving existence of solutions when $F(t,z)$ is bounded on $[-\infty, \max \varphi] \times \Omega$ and continuous in $t$ for each fixed $z \in \Omega$. Later, in \cite{K00}, Kołodziej further generalized this result by considering a function \(F(t,z)\) which is bounded, continuous and non-decreasing in \(t\), and \(d\mu\)-measurable in \(z\). He proved that if there exists a subsolution $v$ satisfying $(dd^c v)^n = d\mu$ together with an appropriate boundary condition then problem \eqref{bt1.0} again has a unique solution. Besides the above mentioned results the Monge-Amp\`ere type equations in Cegrell's classes with given boundary values $\mathcal{F}^a(\O,f)$ or $\mathcal{N}(\O,f)$ have attracted the attention of many authors.  Readers can find results related to this topic through the articles: \cite{Bel14}, \cite{CK06}, \cite{Cz09}. Recently, in the Hermite setting, the Dirichlet problem for Monge-Amp\`ere type equation on compact Hermitian manifolds with boundary has been considered by Ko{\l}odziej and Cuong in \cite{KN23b}. Let $(M,\om)$ be a smooth compact $n$-dimensional Hermitian manifold with the non-empty boundary $\pa M$ and $\mu$ be a positive Radon measure on $M = \overline{M}\setminus\pa{M}$, $\om$ be a Hermitian metric on $M$. Suppose that $F:\mathbb{R}\times M\longrightarrow[0,+\infty)$  is a non-negative function and $\varphi\in C(\pa{M})$. The Dirichlet problem for Monge-Amp\`ere type equations on a Hermitian manifold $M$ with the boundary $\pa M$ in the class of $\om$-plurisubharmonic functions is to find
\[
\begin{cases}
u \in \mathrm{PSH}(M,\omega) \cap L^{\infty}(\overline{M}), \\
(\om + dd^c u)^n = F(u,z) d\mu, \\
\lim\limits_{M \ni z \to x} u(z) = \varphi(x), \quad \forall x \in \partial M.
\end{cases}
\tag{1.3} \label{bt1.3}
\]

\n Under some suitable assumptions on the function $F(t,z)$ and the measure $\mu$, namely, $F(t,z)$ is a bounded nonnegative function which is continuous and non-decreasing in the first variable and $\mu$-measurable in the second one, $\mu$ is a positive Radon measure which is locally dominated by Monge-Amp\`ere measures of bounded plurisubharmonic functions. Then the authors achieved that the Dirichlet problem for Monge-Amp\`ere type equations (1.3) on Hermitian manifolds $M$ with the boundary $\pa M$ in the class of $\om$-plurisubharmonic functions is solved if and only if there exists a bounded subsolution $\underline{u}\in PSH(M,\om)\cap L^{\infty}(M)$ satisfying $\lim\limits_{M\ni z\to x}\underline{u}(z) = \varphi(x)$ for every $x\in\pa{M}$ and 
$$(\om + dd^c\underline{u})^n \geq F(\underline{u},z)d\mu\   \ \text{on $M$}.$$

\n Inspired by the above mentioned researches, in this paper we study the Dirichlet problem on bounded domains in $\mathbb{C}^n$ within the class of $(m-\om)$-subharmonic functions, recently introduced and investigated in \cite{KN23c}, where $\omega$ is a Hermitian metric on $\mathbb{C}^n$. This class is a generalization of the class of $m$-subharmonic functions introduced and investigated by Blocki in \cite{Bl05}. For convenience to readers we recall some basic facts concerning to $m$-subharmonic functions. For further details concerning this class of functions, readers are referred to \cite{Bl05}. Let $\O\subset\cn$ be an open subset and $1 \leq m \leq n$ be an integer. Assume that $\beta = dd^c \|z\|^2$ denotes the canonical K\"ahler form on $\mathbb{C}^n$. In \cite{Bl05}, Błocki introduced the notion of $m$-subharmonic functions as a natural generalization of plurisubharmonic functions, and defined the associated complex $m$-Hessian operator by
\[
H_m(u) := (dd^c u)^m \wedge \beta^{n-m}.
\]
This laid the groundwork for developing a potential theory associated with \(m\)-subharmonicity. Unlike the classical psh case, \(m\)-subharmonic functions are not generally invariant under holomorphic changes of variables, and lack geometric characterizations such as mean value properties. These differences necessitate the development of new techniques tailored to the \(m\)-subharmonic context.

\n Similarly to the case of plurisubharmonic functions, there is also interest in studying the Dirichlet problem on bounded open subsets of $\mathbb{C}^n$ for $m$-subharmonic functions. Li in \cite{Li04} established solvability of the Dirichlet problem for complex $m$-Hessian equations:
\[
\begin{cases}
u \in SH_m(\Omega) \cap L^{\infty}(\overline{\Omega}), \\
(dd^c u)^m \wedge \beta^{n-m} = f\, dV_{2n}, \\
\lim\limits_{\Omega \ni z \to x} u(z) = \varphi(x), \quad \forall x \in \partial \Omega,
\end{cases}
\tag{1.3} \label{bt1.3}
\]
where \(\varphi\) is smooth on \(\partial \Omega\) and \(f\) is a strictly positive, smooth function in \(\Omega\), $SH_m(\O)$ denotes the set of $m$-subharmonic functions on $\O$. Later, Dinew and Kołodziej in \cite[Theorem 2.10]{DiKo} extended this result to the case where \(f \in L^q(\Omega, dV_{2n})\) for some \(q > \frac{n}{m}\). The degenerate case of this problem, where the right-hand side is a general measure, was previously treated by Błocki \cite{Bl05}.

Building further on these results, N. C. Nguyen \cite{C12} extended Kołodziej's subsolution theorem to the setting of bounded \(m\)-subharmonic functions, establishing a subsolution principle for the \(m\)-Hessian equation in analogy with the classical Monge–Ampère case.

More recently, Kołodziej and N. C. Nguyen \cite{KN23c} initiated the study of \(m\)-subharmonic functions associated with a fixed positive Hermitian \((1,1)\)-form \(\omega\) on \(\mathbb{C}^n\). A function \(u: \Omega \rightarrow [-\infty, +\infty)\) is called \(m\)-\(\omega\)-subharmonic (or \(m\)-\(\omega\)-sh for short) if it is upper semicontinuous, belongs to \(L^1_{\text{loc}}(\Omega, \omega^n)\), and for any collection \(\gamma_1, \dots, \gamma_{m-1} \in \Gamma_m(\omega)\), the current
\[
dd^c u \wedge \gamma_1 \wedge \cdots \wedge \gamma_{m-1} \wedge \omega^{n-m} \geq 0
\]
in the sense of distributions. Here, \(\Gamma_m(\omega)\) denotes the positive cone:
\[
\Gamma_m(\omega) := \left\{ \gamma \in \Lambda^{1,1}_{\mathbb{R}} : \gamma^k \wedge \omega^{n-k} > 0 \text{ for } k = 1, \dots, m \right\}.
\]

\n and $\Lambda^{1,1}_{\mathbb{R}}$ is the set of real $(1,1)$-forms. When \(u \in C^2(\Omega)\), the complex Hessian operator relative to \(\omega\) is defined by
\[
H_m(u) := (dd^c u)^m \wedge \omega^{n-m}.
\]
Following inductive methods developed by Bedford–Taylor \cite{BT1} and Błocki \cite{Bl05}, the Hessian operator $H_m(u)$ can be extended to locally bounded \(m\)-\(\omega\)-sh functions as positive Radon measures. More details about the construction and properties of the operator $H_m(u)$ for locally bounded $m-\om$ functions $u$, we refer readers to Section 3 in \cite{KN23c}.

In that same work, Kołodziej and Nguyen also addressed the following Dirichlet problem:
\[
\begin{cases}
u \in SH_{m,\omega}(\Omega) \cap L^{\infty}(\overline{\Omega}), \\
H_m(u) = F(u,z)\, d\mu, \\
\lim\limits_{\Omega \ni z \to x} u(z) = \varphi(x), \quad \forall x \in \partial \Omega,
\end{cases}
\tag{*}
\]
in the special case where \(F(t,z) \equiv 1\). The present paper continues this line of investigation by studying problem (*) with a general right-hand side \(F(t,z)\), under the following conditions:

\begin{itemize}
	\item[(A)] \(F\) is the pointwise limsup on \(\mathbb{R} \times \Omega\) of a sequence of upper semicontinuous functions;
	\item[(B)] There exists a \(\mu\)-measurable set \(X \subset \Omega\) with \(\mu(X) = 0\) such that \(t \mapsto F(t,z)\) is continuous for all \(z \in \Omega \setminus X\);
	\item[(C)] There exists a function \(G \in L^1_{\mathrm{loc}}(\Omega, \mu)\) such that
	\[
	F(t,z) \leq G(z), \quad \forall (t,z) \in \mathbb{R} \times \Omega.
	\]
\end{itemize}

The technical relevance of these assumptions will be clarified in Proposition \ref{md1} in the next section. It is worth noting that, in nearly all previous works in this area, the function \(F(t,z)\) is assumed to be continuous and monotone in the first variable. Here, we allow for a significant relaxation of these conditions. Through the paper, by $SH_{m,\om}(\O)$ we denote the set of $(m-\om)$-subharmonic function on $\O$, $\om$ is a Hermitian metric on $\cn$. Also by $SH^{-}_{m,\om}(\O)$ we denote the set of negative $(m-\om)$-subharmonic functions on $\O$. In this aspect, our first main result reads as follows:
\begin{theorem} \label{th1.1}
	Assume that the following conditions hold:

\n 
(a) There exists $v\in SH_{m, \om}^-(\Om)\cap L^{\infty}(\Om)$ satisfying
$$\lim\limits_{\Om \ni z\to x}v(z)=0\,\forall x\in\pa \Om \ \text{and}\ G\mu\leq H_m(v);$$ 
\n	
(b) For every $m-$polar subset $E$ of $\Om$ we have $\mu (E \cap \{G=0\})=0.$

Then the problem (*) has a solution. Moreover, if $ t \mapsto F(t,z)$ is non-decreasing for every 
	$z \in \O \setminus Y$ where $Y$ is a Borel set with $Cap_m (Y)=0,$ then such a solution $u$ is unique.
\end{theorem}
\begin{remark}{\upshape 
(i) We do {\it not} assume monotonicity of $t \mapsto F(t,z)$ as in \cite{KN23b} and the hypothesis $F(t,z)$ is a bounded function in \cite{KN23b} is replaced by the hypothesis (C) which is more general.

\n (ii) The assumption (b) is needed to guarantee that $\mu$ puts no mass on $m-$polar sets.
 It is not clear to us how this assumption can be removed.}
\end{remark}

\n In addition to the study of the Dirichlet problem, the stability of solutions is also an important and actively investigated problem. In 2002,  Cegrell and Ko{\l}odziej \cite{CK06} considered the case when measure $d\mu= 	(dd^c v)^n,$ where $v$ is a bounded plurisubharmonic function such that $\lim\limits_{\Omega\ni z\to x}v(z)=\varphi (x), \forall x\in\partial \Omega, \int_{\Omega}(dd^c v)^n<\infty$ and $0\leq f_j\leq 1$ is a sequence of $d\mu$-measurable functions satisfying $f_jd\mu$ converge weakly to $fd\mu$ in the sense of measures. Assume that $u_j$ are solutions of Dirichlet problem
	\begin{equation}
	\begin{cases}
		u_j\in PSH(\Om)\cap L^{\infty}(\Om)\\
		(dd^c u_j)^n=f_jd\mu\\
		\lim\limits_{\Omega\ni z\to x}u_j(z)=\va (x), \forall x\in\pa \Om. 
	\end{cases}
\end{equation}
Then,  Cegrell and  Ko{\l}odziej proved that $u_j$ converges in capacity to function $u\in PSH(\Om)\cap L^{\infty}(\Om)$ such that $(dd^c u)^n=fd\mu$ and $\lim\limits_{\Omega\ni z\to x}u(z)=\va (x), \forall x\in\pa \Om.$ Very recently,  Kołodziej and Ngoc Cuong Nguyen in \cite{KN23b} (see Theorem 2.7) extended the result in \cite{CK06} from plurisubharmonic functions to $\o$- plurisubharmonic functions on a compact Hermitian manifold with boundary $(M,\omega)$. However, in \cite{KN23b}, the authors were only able to prove
 stability of solutions in  $L^1$-topology. In general, we known that the convergence in $L^1$ does not imply the continuity of the Hessian operators. In the following theorem, we will prove stability of solution in $m$- capacity. This convergence implies the continuity of corresponding Hessian operators.
\begin{theorem} \label{th1.2}
%Let $\mu, v$ be given as in \ref{th1.1}. 
Let $F, F_j: \R \times \Om \to [0, \infty) \ (j \ge 1)$ be a set
of $dt\times d\mu$- measurable functions that satisfy the conditions (A), (B) and (C) and
that for every $z \in \Om \setminus X,$ the sequence
$F_j (t, z)$ converges locally uniformly to $F (t,z).$
Moreover, suppose that the conditions (a) and (b) in Theorem \ref{th1.1} also hold.
For each $j,$ let $u_j$ be a solution of the equation
	\begin{equation}
\begin{cases}
u_j\in SH_{m, \om}(\Om)\cap L^{\infty}(\Om)\\
H_{m}(u_j)=F_j(u_j,z)d\mu\\
\lim\limits_{\O\ni z\to x}u_j(z)=\va (x), \forall x\in\pa \Om. 
\end{cases}
\end{equation}
Then there exists a subsequence of $\{u_j\}$ that converges in $m$-capacity to a solution $u$ of the problem (*).
\end{theorem}
\n Our work can be viewed as a natural continuation of the recent results  in \cite{AHDP},  obtained for Hessian type equations in the flat background setting. The main contribution of the present paper is to extend the existence and stability theory to the much broader class of $m-\omega$-subharmonic functions, where $\omega$ is an arbitrary Hermitian metric. This extension is far from being merely formal, since in the Hermitian setting the metric is no longer closed. Consequently, one must deal with delicate torsion terms such as $dd^c\omega$ and $d\omega\wedge d^c\omega$, which do not arise in the K\"ahler or flat settings. These additional terms substantially complicate the derivation of a priori estimates, comparison principles, and convergence arguments for Hessian measures. Overcoming these difficulties constitutes one of the main technical achievements of the paper.\\
\n The paper is organized as follows. Besides the introduction, the paper has other two sections. In Section 2, following seminal work \cite{KN23c}, \cite{GN18},... we collect basic features of $m-\om$-sh functions on bounded domains of $\mathbb C^n.$
Most notably is the comparison principle for the Hessian operator $H_m (u)$. Moreover,  we also recall the subsolution theorem, a powerful tool to check existence of solution of Dirichlet problem for Hessian operator. Another important ingredient is the convergence in $m$-capacity of Hessian operator. In Section 3, we supply in details the proofs of our main results.

\n 
{\bf Acknowledgments} 
We wish to express our gratitude to  anonymous referees for their careful reading and constructive comments that help to improve significantly our exposition.\\
Hoang Thieu Anh was funded by the PhD Scholarship Programme of Vingroup Innovation Foundation (VINIF), VinUniversity, code VINIF.2025.TS09.\\
This work is written in our visit in Vietnam Institute for Advanced Study in Mathematics (VIASM) in the Spring of 2025. We also thank VIASM for hospitality.	
	\section{Preliminaries} 
It is important to observe that, by a result of Michelsohn (see equation (4.8) in \cite{Mi82}),
 for $\ga_1,...,\ga_{m-1} \in \Ga_m(\om)$ there is a unique $(1,1)$- positive form $\al$ such that 
$$
\al^{n-1} = \ga_1 \wed \cdots \wed \ga_{m-1} \wed \om^{n-m}.
$$
Hence, the notion of $m-\om$-subharmonicity can be translated, in terms of potential theory, using  the notion of  $\al$-subharmonicity (see e.g., \cite[Definition~2.1, Lemma~9.10]{GN18}). Using this approach, several potential-theoretic properties of $m-\om$-sh functions can be derived from those of $\al$-sh ones. Following \cite{KN23c}, we include below some basic properties of $m-\om$-sh functions. 
\begin{proposition} \label{prop:closure-max}
	Let $\Omega$ be a bounded open set in $\bC^n$. 
	\begin{enumerate}
		\item[(a)]
		If $u_1 \geq u_2 \geq  \cdots$ is a decreasing sequence of $m-\om$-sh functions, then $u := \lim_{j\to \infty} u_j$ is either $m-\om$-sh or $\equiv -\infty$.
		\item[(b)] 
		If $u, v$ belong to $SH_{m, \om}(\Omega)$, then so does $\max\{u,v\}$.
		\item[(c)] 
		If $u, v$ belong to $SH_{m, \om}(\Omega)$ and satisfies $u \le v$ a.e. with respect to Lebesgue measure then $u \le v$ on $\Om.$
		
		\item[(d)] (Theorem 8.7 in \cite{KN23c})
		Let $\{u_\alpha\}_{\alpha \in I} \subset SH_{m,\om}(\Omega)$ be a family locally uniformly bounded above. Put $u(z) := \sup_\alpha u_\alpha(z)$. Then, the upper semicontinuous regularization $u^*$ is $m-\om$-sh. Moreover, the negligible set $E:=\{u<u^*\}$ is $m-$polar, i.e., there exists $u \in SH_{m,\om} (\Om)$ such that $u \equiv -\infty$ on $E.$
	\end{enumerate}
\end{proposition}
\n
Notice that $(c)$ follows from Lemma 9.6 in \cite{GN18} where $u$ and $v$ are viewed as $\al-$subharmonic function with any $(1,1)$-form $\al$ such that 
$
\al^{n-1} = \ga^{m-1} \wed \om^{n-m},$ where $\ga$ is a certain form belonging to $\Ga_m (\om)$.

\n We next recall the notion of $m$-capacity associated with Hessian operators of bounded $m-\o-$subharmonic functions, see Section 4 in \cite{KN23c}.
	\begin{definition}
		For a Borel set $E\subset \Om,$ we set
		$$Cap_m(E)=Cap_m(E,\O):=\sup\left\{\int_{E}(dd^cw)^m \wedge \o^{n-m}: w\in SH_{m,\om}(\Om), 0\leq w\leq 1\right\}.$$
	\end{definition}
\n We recall the following definition as in \cite{KN23c}(Definition {\bf 4.2}).
\begin{definition}
	A sequence of Borel functions $u_j$ in $\Omega$ is said to converge in $m$-capacity (or in $Cap_m(.))$ to $u$ if for any $\delta>0$ and $K\Subset \Omega$ we have 
	$$\lim\limits_{j\to\infty}Cap_m(K\cap|u_j-u|\geq\delta)=0.$$
	\end{definition}
\n We recall Corollary 4.11 in \cite{KN23c} which said that the monotone convergence of locally uniformly bounded sequences of
$m-\om$-sh functions implies convergence in $m$-capacity.
\begin{proposition}\label{conmono}
	Let $\{u_j\}_{j\geq 1}$ be a uniformly bounded and monotone sequence of $m-\om$-sh functions that either $u_j\searrow u$ pointwise or $u_j\nearrow u$ almost everywhere for a bounded  $m-\om$-sh function $u$ in $\Omega.$ Then, $u_j$ converges to $u$ in $m$-capacity.
\end{proposition}
\n
One use of $m$-capacity is to characterize $m-$polar sets.
Recall that, according to Section 7 in \cite{KN23c}, a subset $E$ of $\mathbb C^n$ is said to be $m-$polar if for each $z\in E$ there is an open set $U$ containing $z$ and an $m-\om$-sh function $u$ in $U$ such that $E \cap U \subset \{u=-\infty\}$.
Then by Proposition 7.7 (c) in \cite{KN23c} we know that a subset $E$ of a bounded strictly pseudoconvex domain 
$\Om \subset \mathbb C^n$ is $m$-polar if and only if 
$$Cap^*_m (E):= \inf \{C_m (U): E \subset U, U \subset \Om\ \text{is open}\}=0.$$
\n

\n A major tool in pluripotential theory is the comparison principle.
Before recalling a version of this result for bounded $m-\om-$sh. functions, following Section 6 in \cite{KN23c}, let us fix  a constant ${\bf B}>0$ such that on $\overline \Om$,
$$
- \bb \om^2 \leq dd^c \om \leq \bb \om^2, \quad
-\bb \om^3 \leq d\om \wed d^c \om \leq \bb \om^3.
$$
We also denote by $\rho$ a negative strictly psh. function on a neighborhood of $\overline \Om$ that satisfies $dd^c \rho \geq \om$ on $\Om$.
Now we are able to formulate this important principle (see Theorem 6.1 in \cite{KN23c})
\begin{theorem} \label{lem:wcp-c} 
Let $u,v$ be bounded $m-\om$-sh functions in $\Om$ such that 
$$d = \sup_{\Om} (v-u) >0 \ \text{and}\ \liminf_{z\to \pa \Om} (u-v)(z) \geq 0.$$
Fix $0< \veps < \min\{\frac{1}{2}, \frac{d}{2 \|\rho\|_\infty}\}$. Let us denote for 
$0< s< \veps_0:=\fr{16\veps^3}\bb,$
$$
U(\veps,s):=\{u<(v+\veps\rho) + S(\veps) +s \}, \quad\text{where }
S(\veps)= \inf_\Om[ u - (v+\veps\rho)].
$$
Then,
$$
\int_{U(\veps,s)} H_m (v+\veps\rho)  \leq 
\left( 1+ \frac{C s}{\veps^m} \right) \int_{U(\veps,s)} H_m(u), 
$$
where $C$ is a uniform constant depending on $m,n,\om$.	

\end{theorem}
\n The above theorem has an useful corollary as follow (see Corollary 6.2 in \cite{KN23c}).
\begin{theorem}\label{comparison}
	Let $u,v$ be bounded $m-\om$-sh functions in a neighborhood of $\overline{\Omega}$ such that $\liminf\limits_{z\to\partial\Omega}(u-v)(z)\geq 0.$ Assume that $H_m(v)\geq H_m(u)$ in $\Omega.$ Then  $u\geq v$ on $\Omega.$
	\end{theorem}
\n 

The following result is similar to Corollary 8.4 in \cite{KN23c}. For convenience to readers, we provide the proof here.

\begin{lemma}\label{bd1}
	Assume that $\mu$ vanishes on $m-$polar sets of $\Om$ and $\mu(\Omega) < \infty.$ Let $\{u_{j}\}\in SH_{m, \om}^{-} (\Omega)$ be a sequence satisfying the following conditions:
	
	\n 
	(i) $\sup\limits_{j \ge 1} \int\limits_{\Om} -u_jd\mu <\infty;$
	
	\n 
	(ii) $u_j \to u \in SH_{m, \om}^{-} (\Om)$ a.e. $dV_{2n}.$
	
	Then we have 
	$$\lim_{j \to \infty} \int\limits_{\Om} \vert u_j- u \vert d\mu=0.$$
	In particular $u_j \to u$ a.e. $d\mu$ on $\Om.$
\end{lemma}	
\begin{proof}
	We split the proof into two steps.
	
	\n 
	{\it Step 1.} Firstly, we will prove
	\begin{equation} \label{eq3}
		\lim\limits_{j\to \infty}\int_{\Omega}u_{j}d\mu=\int_{\Omega}ud\mu.
	\end{equation}
	Indeed, in view of (i), by passing to a subsequence we may obtain that 
	\begin{equation} \label{eq7} \lim\limits_{j\to \infty}\int_{\Omega}u_{j}d\mu=a.
	\end{equation}
	By Monotone Convergence Theorem, we have 
	$$\lim_{N \to \infty} \int_{\Om} \max \{u, -N\} d\mu=\int_{\Om} ud\mu,$$	
	and for each $N \ge 1$ fixed 
	$$\lim_{j \to \infty} \int_{\Om} \max \{u_j, -N\} d\mu=\int_{\Om} \max \{u, -N\} d\mu.$$
	Thus, using a diagonal process, it remains to prove (\ref{eq3})	under the restriction that $u_j$ and $u$ are all uniformly bounded from below.
	Since $\mu(\Om)<\infty$  we infer that the set $A:=\{u_j\}_{j \ge 1}$ is bounded in
	the Hilbert space $L^2 (\Om,\mu)$. Therefore, according to Mazur's theorem, we can find a sequence $\tilde u_j$ belonging to the convex hull of $A$ that converges to some element $\tilde u \in L^2 (\Om, \mu).$ After switching to a subsequence we may assume that $\tilde u_j \to \tilde u$ a.e. in $d\mu.$\\
	 On the other hand, it follows from assumption (ii) that $\tilde u_j \to u$ in $L^2 (\Om, dV_{2n})$. This implies that
	$(\sup\limits_{k \ge j} \tilde u_k)^* \downarrow u \ \text{entirely on}\  \Om.$ Moreover, according to Theorem 7.8 of \cite{KN23c}, we see that $m-$negligible set $\{(\sup\limits_{k \ge j} \tilde u_k)<(\sup\limits_{k \ge j} \tilde u_k)^*\}$ are $m-$polar set.
	Therefore, using Monotone Convergence Theorem we get
	$$
	\int_{\Om} ud\mu=\lim_{j \to \infty} \int_{\Om} (\sup\limits_{k \ge j} \tilde u_k)^* d\mu
	=\lim_{j \to \infty}\int_{\Om} (\sup\limits_{k \ge j} \tilde u_k) d\mu
	=\int_{\Om} \tilde ud\mu =a.$$
	Here the second equality follows from the fact that $\mu$ does not charge $m-$polar  sets  
 and the last equality results are obtained from the choice of $\tilde u_j$ and (\ref{eq7}). The equation (\ref{eq3}) follows.
	
	\n 
	{\it Step 2.} Completion of the proof. We put $v_j:= (\sup\limits_{k \ge j} u_k)^*$. Then we have $v_j \ge u_j, v_j \downarrow u$ on $\Omega$ and $v_j \to u$ in $L^1 (\Om, dV_{2n}).$
	So by the result obtained in Step 1 we have
	\begin{equation} \label{eq4}
		\lim_{j \to \infty} \int_{\Om} v_j d\mu= \int_\Om ud\mu=\lim_{j \to \infty} \int_{\Om} u_j d\mu. 
	\end{equation} 
	Using the triangle in equality we obtain
	$$\begin{aligned} 
		\int_{\Om} \vert u_j-u\vert d\mu &\le \int_{\Om} (v_j-u)d\mu+ \int_{\Om} (v_j-u_j)d\mu\\
		&=  2\int_{\Om} (v_j-u)d\mu + \int_{\Om} (u-u_j)d\mu.
	\end{aligned}$$
	Hence by applying (\ref{eq4}) we finish the proof of the lemma.
\end{proof}	
\n According to Lemma 5.1 and Lemma 5.4 in \cite{KN23c}, we know that $H_m (u)$ is continuous
with respect to  monotone convergent of locally uniformly bounded sequences in $SH_{m,\om} (\Om).$
In our work, this fact will be referred to as the monotone convergence theorem.
We present below an analogous result where monotone convergent is replaced by convergence in $m$-capacity.
Of course, this fact is inspired by a well known result of Xing in \cite{Xi96} as in the classical case of plurisubharmonic functions.
\begin{proposition} \label{cap}
Let $u_j$ be a locally uniformly bounded sequence in $SH_{m,\om} (\Om).$ Assume that $u_j$ converges in $m$-capacity
to  a locally bounded $u \in SH_{m,\om} (\Om)$. Then $H_m (u_j)$ converges weakly to $H_m (u).$
\end{proposition}
\begin{proof}
\ Because the problem is local we may assume that $\Omega=\mathbb{B}$ is a ball in $\cn$. Moreover, we assume $-1\leq u_j, u\leq 0$ on $\mathbb{B}$ for all $j\geq 1$. Let $\va$ be a test function on $\mathbb{B}$. By the localization principle (see Section 2 in \cite{KN23c}) there exists a fixed compact subset $A\Subset\mathbb{B}$ such that $supp\va\subset A$ and $u_j=u$ on $\mathbb{B}\setminus A$. We have to show
$$\lim_{j \to \infty} \int \va H_m (u_j) = \int \va H_m (u).$$ 
\n We have 
$$H_m (u_j)-H_m (u)= dd^c (u_j-u)\wedge \sum_{s=0}^{m-1} (dd^c u_j)^s \wedge (dd^c u)^{m-1-s} \wedge \om^{n-m}.$$

\n Thus it follows that 
\begin{align*}
	&\int\limits_{\mathbb{B}}\va\Bigl[H_m(u_j)- H_m(u)\Bigl]=\int\limits_{\mathbb{B}}\va dd^c(u_j-u)\wedge\sum\limits_{s=0}^{m-1}(dd^c u_j)^s\wedge(dd^c u)^{m-1-s}\wedge\om^{n-m}\\
	&=\int\limits_{\mathbb{B}}\va \om^{n-m}dd^c(u_j-u)\wedge\sum\limits_{s=0}^{m-1}(dd^c u_j)^s\wedge(dd^c u)^{m-1-s}\hskip4cm(2.1)
\end{align*}

\n By Stoke's formula we infer that the right-hand side of (2.1) is equal 
$$(2.1)= \int\limits_{\mathbb{B}}(u_j-u)dd^c(\va\om^{n-m}) \wedge \sum\limits_{s=0}^{m-1}(dd^c u_j)^s\wedge(dd^c u)^{m-1-s}.\eqno(2.2)$$

\n On the other hand, by Corollary 2.4 in \cite{KN16} we have
\begin{align*}
	&\Bigl|dd^c(\va\om^{n-m}) \wedge \sum\limits_{s=0}^{m-1}(dd^c u_j)^s\wedge(dd^c u)^{m-1-s}\Bigl|\\
	&\leq C\Bigl[dd^c (su_j+ (m-1-s)u)\Bigl]^{m}\wedge\om^{n-m}.\hskip5cm(2.3)
\end{align*}

\n where $C>0$ is a constant which is dependent on $\om, \va$ and is independent of $j.$  Given $\varepsilon>0$. Coupling (2.1) and (2.2) we infer that
\begin{align*}
	&\int\limits_{\mathbb{B}}\va\Bigl[H_m(u_j)- H_m(u)\Bigl]=\int\limits_{\mathbb{B}}(u_j-u)dd^c(\va\om^{n-m}) \wedge \sum\limits_{s=0}^{m-1}(dd^c u_j)^s\wedge(dd^c u)^{m-1-s}\\
	&=\int\limits_{A\cap\{|u_j-u|< \varepsilon\}}(u_j-u)dd^c(\va\om^{n-m})\wedge\sum\limits_{s=0}^{m-1}(dd^c u_j)^s\wedge(dd^c u)^{m-1-s}\\
	&+\int\limits_{A\cap\{|u_j-u|\geq \varepsilon\}}(u_j-u)dd^c(\va\om^{n-m})\wedge\sum\limits_{s=0}^{m-1}(dd^c u_j)^s\wedge(dd^c u)^{m-1-s}. 
\end{align*}

\n However, $-m< su_j +(m-1-s)u\leq 0$. Hence, by (2.3) we get that
\begin{align*}
	&\Biggl|\int\limits_{\mathbb{B}}\va\Bigl[H_m(u_j)- H_m(u)\Bigl]\Biggl|\leq\\
	&\leq \int\limits_{A\cap\{|u_j-u|< \varepsilon\}}|(u_j-u)||dd^c(\va\om^{n-m})\wedge\sum\limits_{s=0}^{m-1}(dd^c u_j)^s\wedge(dd^c u)^{m-1-s}|+\\
	&+\int\limits_{A\cap\{|u_j-u|\geq \varepsilon\}}|(u_j-u)||dd^c(\va\om^{n-m})\wedge\sum\limits_{s=0}^{m-1}(dd^c u_j)^s\wedge(dd^c u)^{m-1-s}|\\
	&\leq Cm^{m}Cap_{m}(A)\varepsilon + 2Cm^{m}Cap_{m}\Bigl(A\cap\{|u_j-u|\geq\varepsilon\}\Bigl).
\end{align*}
\n By the hypothesis, $Cap_{m}\Bigl(A\cap\{|u_j-u|\geq\varepsilon\}\Bigl)\longrightarrow 0$ as $j\to\infty$. Therefore, we get that $H_m(u_j)$ is weak*-convergent to $H_m(u)$. The proof is complete.
\end{proof}
\n To see that the problem  (*) is well posed, we need the following elementary fact for which no originality is claimed.
\begin{proposition}\label{md1}
	Let
	$F(t,z): \R\times \Om \to [0,+\infty)$ be a function
	that satisfies $(A), (B), (C)$ and $u$ be a locally bounded upper semicontinuous function on $\O.$ Then 
	$F(u,z)d\mu$ is a positive Radon measure on $\O.$
\end{proposition}
\begin{proof}
	First, we show that $F(u,z)$ is $\mu-$ measurable on $\O \setminus X$.	
	Let $u_k$ be a sequence of continuous functions on $\O$ that decreases to $u$ pointwise on $\O.$
	It follows that $F(u_k, z)$ converges pointwise to $F(u,z)$ for all $z \in \O \setminus X.$
	Thus, it suffices to check that $F(u_k,z)$ is $\mu-$measurable.
	Since, by the assumption, $F$ is the pointwise limsup of a sequence of upper semicontinuous functions on $\mathbb R \times \O,$ we may assume $F$ is upper semicontinuous. We claim that $F(u_k, z)$ is upper semicontinuous on $\O.$ Indeed, fix a sequence $\O \ni \{z_j\} \to z^* \in \O.$
	Then $u_k (z_j) \to u_k (z^*)$ as $j \to \infty.$ Hence $(u_k (z_j), z_j) \to (u_k (z^*), z^*)$ as $j \to \infty.$
	Therefore
$$\limsup_{j \to \infty} F(u_k (z_j), z_j) \le F(u_k (z^*), z^*).$$
Hence $F(u_k,z)$ is non-negative upper semicontinuous on $\O$ as claimed.
It follows that $F(u,z)$ is $\mu-$ measurable on $\O \setminus X.$ 
To finish off, we let $\va \ge 0$ be a continuous function with compact support in $\O.$
Since $u$ is bounded on $\O$, we have
$$0 \le M:= \sup \{|u_k (z)|, |u(z)|: z \in K:= \text{supp} \va, k \ge 1\}< \infty.$$
By the assumption (C), we have
$$F(t, z) \le G(z) \ \forall (t,z) \in [-M,M] \times K \Rightarrow F(u_k (z),z) \le G(z) \ \forall z \in K.$$
Thus, using	Lebesgue dominated convergence theorem we conclude that 
	$$0 \le \lim_{k \to \infty} \int\limits_{\O \setminus X} \va (z) F(u_k, z)d\mu=\int\limits_{\O \setminus X} \va (z) F(u,z)d\mu=\int\limits_{\O} \va (z) F(u,z)d\mu$$
 Therefore, by the Riesz representation theorem, $F(u,z)d\mu$ can be identified with a positive Radon measure on $\Om$ as required.
\end{proof}	
\begin{remark}
We do not know if the condition (A)	can be relaxed to $F$ is just a $dt \times d\mu$ measurable function as 
mostly assumed in the literature.
\end{remark}
\n
Now we formulate the following subsolution theorem  which plays a prominent role in our work (see Theorem 8.7 in \cite{KN23c}).
\begin{theorem}\label{sub}
Assume that there exists $v \in SH_{m,\om} (\Om) \cap L^\infty (\Om)$ satisfying 
\[\label{eq:bounded-subsol}
H_m(v) \geq \mu, \quad \lim_{x\to \d \Om}  v(x) = 0.
\]	
Then, there exists a unique bounded $m-\om$-sh function $u$ solving  
$$\lim_{z\to x} u(z) = \va (x), \forall  x\in\d\Om,
	H_m(u) = \mu\ \text{on } \Om.
	$$
\end{theorem}
	\section{ Weak solution to  Hessian type equation}
We need the following version of the comparision principle. A similar result was obtained 
for quasi-plurisubharmonic  functions in Proposition 2.2 of  \cite{KN23b}.
	\begin{proposition}\label{md4}
Let $\nu \ge \mu$ be positive Radon measures on $\Om.$ 
Assume that	$ t \mapsto F(t,z)$ is a non-decreasing function in $t$ 
for all $z \in \Om \setminus Y,$ where $Y \subset \Om$ is a Borel set with $C_m (Y)=0.$
Let
		$u,v\in SH_{m, \om}(\Om)\cap L^{\infty}(\Om)$ be functions satisfying the following conditions: 

\n 
(i)
$\liminf\limits_{z\to\pa \Om}(u-v)(z)\geq 0;$

\n
(ii) $H_{m}(u)=F(u,z)\mu, H_{m}(v)=\tilde F(v,z)\nu,$ where $\tilde F \ge F$ is a
measurable function on $\Om.$
  
 Then $u\geq v$ on $\Om.$
	\end{proposition}
	\begin{proof}
		 By subtracting from $u,v$ a constant $c$ and replacing $F(t,z)$ and $\tilde F(t,z)$ by $F(t+c,z)$ and $\tilde F(t+c,z)$ we may also assume that $u,v< 0.$
		Arguing by contradiction, suppose that $\{u<v\}\neq \emptyset$. Now we will proceed in the same way as Corollary 6.2 in \cite{KN23c}. For readers convenience, we give some details.
Set $$d:=\sup_{\Om}(v-u)>0.$$
Pick $z_0 \in \Om$ such that
$$v(z_0)-u(z_0)>\fr{2d}3.$$
Hence, there exist positive constants $ \de, a, b$ so small that:

\n 
(a) $\de>a\vv v\vv_{\infty};$

\n 
(b) $\fr{d}6>a \vv v\vv_{\infty}+\de;$

\n 
(c) $\fr{d}6 \ge a \vv v\vv_{\infty} +b \vv \rho\vv_{\infty}.$

\n
By (a) and the assumption (i), we infer that
 $$\liminf_{z\to \partial \Om} [(u+\de) - (1+ a)v](z) \geq \de -a \vv v\vv_{\infty}>0.$$
Next, from (b) we get
 $$\begin{aligned}
 \sup\limits_{\Om} [(1+ a) v - (u+\de)] &\ge (1+a)v(z_0)-(u(z_0)+\de)\\
 &=(v(z_0)-u(z_0))+a v(z_0)-\de\\
 & \ge \fr{2d}3 - (a \vv v\vv_{\infty}+\de)\\
 &>\fr{d}2.
 \end{aligned}$$
Fix 
$$0< \veps < \min \Big \{\frac{1}{2}, \frac{d}{4 \|\rho\|_\infty}, b \Big \}, 0< s< \min \Big \{\fr{d}3, \veps_0:=\fr{16\veps^3}\bb \Big \}.$$
\n
In view of the above estimates, we may apply Theorem~\ref{lem:wcp-c} to the functions
$$\tilde u:= u+\de, \tilde v: = (1+a)v$$
to obtain for $0<s<\veps_0$, the following estimate
		\begin{equation} \label{comp}
		\int_{\tilde U(\veps,s)} H_m(\tilde v+\veps\rho) \leq \left( 1+ \frac{C s}{\veps^m} \right) \int_{\tilde U(\veps,s)} H_m(u),
		\end{equation}
where		
$$
\tilde U(\veps,s):=\{\tilde u<(\tilde v+\veps\rho) + S(\veps) +s \}, 
S(\veps):= \inf_\Om[ \tilde u - (\tilde v+\veps\rho)].$$
Notice that, being a non-empty set, $\tilde U(\veps,s)$, has a positive Lebesgue measure, according to Lemma 9.6 in \cite{GN18}. 	
Observe that
$$\begin{aligned} 
S(\ve)&=\inf_{z \in \Om}	[u(z)+\de-(1+a)v(z)-\ve \rho (z)]\\
&\le  u(z_0)+\de-(1+a) v(z_0)-\ve \rho(z_0)\\
&\le \de+v(z_0)-\fr{2d}3-(1+a)v(z_0)+\ve \vv \rho \vv_{\infty}\\
&\le\de-\fr{2d}3+a\vv v\vv_{\infty}+\ve \vv \rho \vv_{\infty}\\
&\le \de-\fr{d}2,
\end{aligned}$$
where the last estimate results from (c).
Now we claim that 
\begin{equation} \label{eqi}
\tilde  U(\ve,s) \subset \{z \in \Om: u(z)<v(z)\}.
\end{equation}
Assume otherwise, then there exists $z_1 \in \tilde U(\ve,s)$ but $u(z_1) \ge v(z_1).$\\
It follows that 
$$(1+a)v(z_1)+\ve \rho (z_1)+S(\ve)+s>u(z_1)+\de \ge v(z_1)+\de.$$
Hence
$$\de \le a\vv v\vv_{\infty} +\ve \vv \rho \vv_{\infty}+S(\ve)+s\le a\vv v\vv_{\infty} +b \vv \rho \vv_{\infty}+\de-\fr{d}2+s.$$
This yields
$$\fr{d}2 \le  a\vv v\vv_{\infty} +b \vv \rho \vv_{\infty}+s \le \fr{d}6+s.$$
where the last estimate follows from (c). We thus obtain a contradition to the choice of $s.$ 
Hence we have proved the inclusion (\ref{eqi}). Therefore, on  $\tilde U(\ve,s) \setminus Y,$ using (ii) we obtain
$$H_{m}(u)=F(u,z)d\mu\leq F(v,z)d\mu\leq \tilde F(v,z)d\nu=H_{m}(v).$$
Since $C_m (Y)=0,$ and since $H_m (u)$ and $H_m (v)$ does not charge $Y$, we see that
$H_m (u) \le H_m (v)$ entirely on $\tilde U(v,\ve).$
So we get the following estimate on $\tilde U(\ve,s)$
$$H_m (\tilde v+\ve \rho) \ge (1+a)^m H_m (v)+\ve^m H_m (\rho) \ge (1+a)^mH_m (u)+\ve^m H_m(\rho).$$		
Combining this estimate and (\ref{comp}) we obtain 		
$$\int_{\tilde U(\veps,s)} H_m(\rho) \le 0$$
for $s>0$ so small that $(1+ a)^m \geq 1+ Cs/\veps^m$. 	This forces $\tilde U(\ve,s)$ has Lebesgue measure $0.$ 
We arrive at a contradiction.
	\end{proof}
	Now we will prove Theorem \ref{th1.1}.

	\begin{proof} [{\it The proof of Theorem \ref{th1.1}}]
First, we show that $\mu$ puts no mass on $m-$polar subsets of $\Om.$ Indeed, for every  $m-$polar subset $E$ of $\Om,$
 by the assumption (a) we infer that $(G\mu) (E)=0.$ Combining this result with hypothesis (b), we obtain $\mu(E)=0.$
		
\n Next, by Theorem \ref{sub}, there exists $h\in SH_{m, \om}(\Om)\cap L^{\infty}(\Om)$ such that
\begin{equation} \label{eqh}
H_{m}(h)=0,\, h=\va \ \text{on}\ \pa \Om.
\end{equation}

\n We put
$$A:= \max \{ -\inf_{\Omega} (v+h), 0\}.$$ 
%Then by rescaling, we may assume $0 \le F \le 1$ on $[-A, A] \times \O.$		
Consider the set
$$\mathcal{A}:=\{u\in SH_{m, \om}(\Om):v+h\leq u\leq h\}.$$
Since $h\in\mathcal{A},$ we infer that $\mathcal{A}\neq\emptyset.$
Notice also that 
$$-A \le u(z) \le A, \ \forall z \in \Omega, \forall u \in \mathcal A.$$	
It is also easy to see that $\mathcal{A}$ is a convex and bounded set in $L^1(\Om)$ with respect to $L^1-$topology. Hence, $\mathcal A$ is a convex compact set in $L^1 (\Om).$

Moreover, accordingt to  Proposition \ref{md1}, for $u\in\mathcal{A},$ we have
$F(u,z)d\mu$ is a positive Radon measure. By condition (C), we get $F(u,z)d\mu\leq Gd\mu\leq H_m(v).$ Hence, by Theorem \ref{sub}, we obtain a unique function  $g\in SH_{m, \om}(\Om)\cap L^{\infty}(\Om)$ satisfying
$$
H_{m}(g)=F(u,z)d\mu, \lim_{z\to x} g(z) = \va (x), \forall  x\in\d\Om$$ 
Obviously, we have
\begin{equation}\label{eq6.3}H_{m} (v+h) \ge H_{m} (v) \ge G.d\mu \ge F(u,z)d\mu=H_{m} (g)\geq H_{m}(h)=0.\end{equation}
According to Theorem \ref{comparison}, we obtain $h\geq g \ge v+h$ on $\Om.$ Thus, we deduce that $g\in\mathcal{A}$.	
Therefore, we can define the map $$T:\mathcal{A}\to\mathcal{A}, T(u):= g.$$

\n Next, we will verify that $T$ is continuous. Indeed, let $\{u_j\}\subset \mathcal{A}$ be a sequence such that $u_j\to u$ in $L^1(\Om).$
% We put $g=T(u)$ and $g_j=T(u_j).$ 
Since $v+h\leq u,u_j\leq h$, by Lemma \ref{bd1}, we obtain $$\lim\limits_{j}\int_{\Omega} |u_j-u|d\mu=0.$$ 
Hence $u_j \to u$ almost everywhere ($d\mu$).		
Now we define for $z \in \Omega,$ the following sequences of non-negative uniformly bounded measurable functions
$$\theta^1_j (z):= \inf_{k \ge j} F(u_k (z),z),
\theta^2_j (z):= \sup_{k \ge j} F(u_k (z),z).$$
Since $F(t,z)$ is continuous function in the first variable, we have:\\
\n 
(i) $0 \le \theta^1_j (z) \le  F(u_j (z),z) \le \theta^2_j (z) \le G$ for $j \ge 1;$

\n 
(ii) $\lim\limits_{j \to \infty} \theta^1_j (z)=\lim\limits_{j \to \infty} \theta^2_j (z)= F(u(z),z)$ almost everywhere corresponding to the measure ($d\mu$).

\n It follows from Theorem \ref{sub} that we can find $\gamma^1_j, \gamma^2_j \in SH_{m, \om}(\Om)\cap L^{\infty}(\Om)$ are solutions of the equations 
$$H_{m} (\gamma^1_j)=\theta^1_j \mu, \lim_{z\to x} \gamma^1_j(z) = \va (x), \forall  x\in\d\Om$$ 
$$H_{m} (\gamma^2_j)=\theta^2_j \mu, \lim_{z\to x} \gamma^2_j(z) = \va (x), \forall  x\in\d\Om.$$
Using the same argument as in inequality \eqref{eq6.3}, we also have $v+h\leq \gamma^1_j,\gamma^2_j\leq h.$ 
Note that, we have $\{\theta^1_j\}_j$ is a increasing sequence and $\{\theta^2_j\}_j$ is a decreasing sequence. Then, using Theorem \ref{comparison} we see that $\gamma^1_j \downarrow \gamma^1\in  SH_{m, \om}(\Om)\cap L^{\infty}(\Om)$ and $\gamma^2_j \uparrow (\gamma^2)^*$ outside an $m$-polar set with $(\gamma^2)^*\in  SH_{m, \om}(\Om)\cap L^{\infty}(\Om)$. Furthermore, in view of (i) and Proposition \ref{md4} we also have
\begin{equation} \label{e2}
	\gamma^1_j \ge T(u_j) \ge \gamma^2_j.
\end{equation}

\n On the other hand, it follows from (ii) that 
$$H_{m} (\gamma^1_j) \to F (u,z) d\mu$$ and $$ H_{m} \gamma^2_j) \to F(u,z)d\mu.$$

\n Note that since $\gamma^1_j\searrow \gamma^1,$ according to Proposition \ref{conmono},  we have $\gamma^1_j$ converges in $m$-capacity to $\gamma^1.$ By Proposition \ref{cap}, we obtain $H_{m}(\gamma^1_j)$ converges weakly to $H_{m}(\gamma^1).$ Similarly,  we also have $H_{m}(\gamma^2_j)$ converges weakly to $H_{m}(\gamma^2)^*.$ Therefore,  we infer that
$$H_{m} (\gamma^1)=H_{m} ((\gamma^2)^*)= F(u(z),z)d\mu=H_{m}(T(u)).$$
Applying again Theorem \ref{comparison}, we infer that $\gamma^1=(\gamma^2)^*=T(u)$ on $\O.$ 
By the squeezing property (\ref{e2}), we infer that $T(u_j) \to T(u)$ pointwise outside a $m$-polar set of $\Omega.$
Since $\mu$ puts no mass on $m$-polar sets, we may apply Lebesgue dominated convergence theorem
to achieve that $T(u_j) \to T(u)$  in $L^1 (\Omega, d\mu)$. Thus $T: \mathcal A \to \mathcal A$ is continuous. 	
According to Schauder fixed - point theorem, there exists $\tilde{u}\in\mathcal{A}$ such that $T(\tilde{u})=\tilde{u}.$
Therefore $\tilde{u}$ is a solution of the problem (*). \\
\n In the case, $ t \mapsto F(t,z)$ is non-decreasing for every 
$z \in \O \setminus Y$ with $Cap_m (Y)=0,$ then the uniqueness of $\tilde{u}$ follows directly from Proposition \ref{md4} (with $F=\tilde F, \mu=\nu$).

\end{proof}	
	\begin{proof} [{\it The proof of Theorem \ref{th1.2}}]
Let $h\in SH_{m, \om}(\Om)\cap L^{\infty}({\Om})$ be a function that satisfies equation (\ref{eqh}). Then by the assumption on $G$, for each $j \ge 1$ we have 
$$H_m (v+h) \ge H_m (v) \ge Gd\mu \ge F(u_j, z)d\mu=H_m (u_j) \ge H_m (h).$$
Thus using the comparision principle we obtain 
$$v+h \le u_j \le h.$$
In particular, the sequence $u_j$ is uniformly bounded in $SH_{m,\om} (\Om) \cap L^\infty (\Om).$ Hence, after passing to a subsequence, we may assume that
$u_j$ converges pointwsise a.e. $(dV_{2n})$ to a function $u \in SH_{m,\om} (\Om).$ By Lemma \ref{bd1}, we obtain $$\lim\limits_{j}\int_{\Omega} |u_j-u|d\mu=0.$$ 
Hence $u_j \to u$ almost everywhere ($d\mu$).\\

\n Observe that, by the assumption that $F_j (t, \cdot) \to F(t, \cdot)$ locally uniformly on $\R$ when $z \in \Om \setminus X$ is fixed, in view of the condition (B) we get that 
$$\lim_{j \to \infty} F_j (u_j (z), z)= F(u(z), z), \text{a.e.}\ d \mu.$$
\n
Now for $z \in \Om,$ we define 
$$\theta^1_j (z):= \inf_{k \ge j} F_k (u_k (z),z),
\theta^2_j (z):= \sup_{k \ge j} F_k (u_k (z),z).$$
Then, by the assumption (C) we have:\\
\n 
(i) $0 \le \theta^1_j (z) \le  F_j (u_j (z),z) \le \theta^2_j (z) \le G(z)$ for $j \ge 1;$

\n 
(ii) $\lim\limits_{j \to \infty} \theta^1_j (z)=\lim\limits_{j \to \infty} \theta^2_j (z)= F(u(z),z)$ a.e. ($d\mu$).

	\n 	
We apply Theorem \ref{sub} to get $\gamma^1_j, \gamma^2_j, \tilde u \in SH_{m,\om} (\Om) \cap L^\infty (\Om)$ which are solutions of the equations 
$$H_{m} (\gamma^1_j)=\theta^1_j d\mu, \lim_{z\to x} \gamma^1_j(z) = \va (x), \forall  x\in\d\Om$$
and	
$$H_{m} (\gamma^2_j)=\theta^2_j d\mu, \lim_{z\to x} \gamma^2_j(z) = \va (x), \forall  x\in\d\Om$$
and
$$H_{m}( \tilde u)=F(u,z)d\mu, \lim_{z\to x} \tilde u(z) = \va (x), \forall  x\in\d\Om.$$ 	
Now, we may repeat the proof of Theorem \ref{th1.1} to see that 
$\gamma^1_j \downarrow \gamma^1\in SH_{m,\om} (\Om) \cap L^\infty (\Om)$ and $ \gamma^2_j \uparrow (\gamma^2)^*$ outside an $m$-polar set and 
$\gamma^1=(\gamma^2)^*= \tilde u$ on $\O.$ 

\n	Moreover, note that $H_{m}(u_j)=F_j(u_j,z)d\mu,$ in the view of (i) and Theorem \ref{comparison} we have \begin{equation}\label{eq6.5}\gamma^1_j\geq u_j\geq \gamma^2_j.\end{equation} This implies that
\begin{equation} \label{eq1}
	|u_{j}- \tilde u| \le \max \{\gamma^1_j-\tilde u, \tilde u-\gamma^2_j\}
\end{equation}	
Now we claim that $u_j \to \tilde u$ in $m$-capacity on $\Om.$ Indeed, fix $\de>0,$ by (\ref{eq1}), for each compact subset $K$ of $\Om$ we have
$$Cap_m (K \cap \{|u_{j} -\tilde u|>\de\}) \le Cap_m(K \cap \{\gamma^1_j-\tilde u>\de\})+
Cap_m (K \cap \{\tilde u- \gamma^2_j >\de\}).$$
According to Proposition \ref{conmono} we see that both monotone sequences $\gamma^1_j$ and $\gamma^2_j$ converge
to $\tilde u$ in $m$-capacity. Therefore
$$\lim_{j \to \infty} Cap_m (K \cap \{|u_{j} -\tilde u|>\de\})=0.$$
Hence $u_j$ tends to $\tilde u$ in $m$-capacity as claimed.\\
On the other hand, by the monotone convergence theorem, we have
$$\lim_{j \to \infty} \int \gamma^1_j d\mu=\lim_{j \to \infty} \int \gamma^2_j d\mu=\int \tilde u d\mu.$$
Notice that, for the second equality, we use the fact that $\mu$ puts no mass on $m-$polar set, so that $\gamma^2_j \uparrow \tilde u$ a.e. ($d\mu$).
It follows from inequality \eqref{eq6.5} that
$$\lim_{j \to \infty}\int |u_{j}-\tilde u|d\mu=0.$$
Hence, $\tilde u=u$ a.e. $(d\mu),$ because $u_j \to u$ a.e. $(d\mu).$
Now, we have $$H_{m}(\tilde u)=F(u,z)d\mu=F(\tilde u,z)d\mu.$$
It means that we have $\tilde u$ is the solution of problem (*) and we have
$u_j \to \tilde u$ in $m$-capacity on $\Om.$
The proof is complete.

\end{proof}		

	\section*{Declarations}
	\subsection*{Ethical Approval}
	This declaration is not applicable.
	\subsection*{Competing interests}
	The authors have no conflicts of interest to declare that are relevant to the content of this article.
	\subsection*{Authors' contributions }
Hoang Thieu Anh,	Le Mau Hai, Nguyen Quang Dieu and Nguyen Van Phu  together studied  the manuscript.
	%\subsection*{Funding }
	%No funding was received for conducting this study.
	\subsection*{Availability of data and materials}
	This declaration is not applicable.


\begin{thebibliography}{000000}
			\bibitem[AHDP]{AHDP}	Hoang Thieu Anh, Le Hai Mau, Nguyen  Dieu Quang and Nguyen Van Phu, {\it Bounded solution to Hessian type equations for $(\omega,m)-\beta$-subharmonic functions on a ball in $\mathbb{C}^n$}, Journal of Mathematical Analysis and Applications, {\bf 554}(2026), issue 2, Article: 129954. 
		
			\bibitem[BT76]{BT76} E. Bedford and B. A.Taylor, {\it  The Dirichlet problem for a complex Monge–Amp\`ere operator}. Invent. Math, {\bf 37}(1976), 1–44. 
		
		\bibitem[BT79]{BT79} E. Bedford and B. A.Taylor, {\it The Dirichlet problem for an equation of complex Monge-Amp\`ere type}, In: Proceedings of the Partial differential equations and geometry. Lecture Notes in Pure and Appl. Math., 48, pp. 39-50, Park City, Utah, Dekker, New York (1979). 
		
		
		
		\bibitem[BT82] {BT1} E. Bedford and B. A.Taylor, {\sl A new capacity for plurisubharmonic functions,} Acta Math. {\bf 149} (1982), 1-40. https://doi.org/10.1007/BF02392348.
		
		\bibitem[Bel14] {Bel14} S. Benelkourchi, {\it  Weak solution to the complex Monge–Amp\`ere equation on hyperconvex domains}, Ann. Polon. Math. {\bf 112(3)}(2014), 239–246
	
		\bibitem[Bl05]{Bl05} Z. B\l ocki, {\it Weak solutions to the complex Hessian equation,} Ann. Inst. Fourier (Grenoble) {\bf 55} (2005), no.~5, 1735-1756.
		
		\bibitem [CePe92]{CePe92} U.Cegrell and L.Persson {\it The Dirichlet problem for the complex Monge–Amp\`ere operator: stability  in $L^2$}. Michigan Math. J, {\bf39}(1992),145–151.
	
		\bibitem[Ce84]{Ce84} U. Cegrell,  {\it On the Dirichlet problem for the complex Monge-Amp\`ere operator}, Math. Z. {\bf185}(1984), 247–251.
		
			\bibitem[CK06]{CK06}
		U. Cegrell, S.  Ko{\l}odziej, {\sl The equation of complex Monge-Amp\`ere type and stability of solutions}, Math. Annalen,  {\bf 334} (2006), no. 4, 713–729. 
		
		\bibitem[Cz09]{Cz09} R. Czyz, The complex Monge–Amp\`ere operator in the Cegrell classes, Dissertationes  Math.{\bf 466}(2009) 83 pp.
		
	\bibitem[C12]{C12}	Ngoc Cuong Nguyen, {\sl Subsolution theorem for the complex Hessian equation}, Universitatis Iagellonicae Acta Mathematica, {\bf 50} (2012), 69–88. 
		
	\bibitem[DK14]{DiKo} S. Dinew and S. Ko{\l}odziej, {\it A priori estimates for the complex Hessian equations}, Analysis \& PDE, {\bf 7} (2014), 227-244.
		
		\bibitem [GN18]{GN18} D. Gu, N.-C. Nguyen, {\it The Dirichlet problem for a complex Hessian equation on compact Hermitian manifolds with boundary}, Ann. Sc. Norm. Super. Pisa Cl. Sci., {\bf 18} (2018), no. 4, 1189-1248.
		
			\bibitem[Ko96]{Ko96} S. Ko{\l}odziej, {\it Some sufficient conditions for solvability of the Dirichlet problem for the complex  Monge–Amp\`ere operator}.Ann. Polon. Math. {\bf65}(1996), 11–21.
			
				\bibitem[Ko98]{Ko98} S. Ko{\l}odziej, {\it The complex  Monge–Amp\`ere equation}. Acta Math. {\bf 80} (1998), 69-117.
		
		\bibitem[K95]{K95} S. Ko{\l}odziej, {\it The range of the complex Monge-Amp\`ere operator.II}, Indiana Univ. Math.,  {\bf 44}, no. 3, 765-782 (1995).
		
			
		\bibitem[K00]{K00} S. Ko{\l}odziej, {\it Weak solutions of equations on complex Monge-Amp\`ere type}, Ann. Polon. Math.,  {\bf 73}, no. 1, 59-67 (2000).
		
	\bibitem[KN16]{KN16} S. Ko{\l}odziej and N.-C. Nguyen, {\it Weak solutions of complex Hessian equations on compact Hermitian manifold}, Compos. Math.,  {\bf 152} (2016), 2221-2248.
		
		\bibitem[KN23a] {KN23a} S. Ko{\l}odziej and N.-C. Nguyen, {\it  The Dirichlet problem for the Monge-Amp\`ere equation on Hermitian manifolds with boundary}, Calc. Var. Partial Differential Equations., {\bf 62} (2023), no. 1, Paper number 1.
		
		\bibitem[KN23b] {KN23b} S. Ko{\l}odziej and N.-C. Nguyen, {\it Weak solutions to Monge-Amp\`ere type equations on compact Hermitian manifold with boundary}, J. Geom. Anal., {\bf 33} (2023), no. 1, Paper number 11, 20 pp.
		
	\bibitem[KN23c] {KN23c} S. Ko{\l}odziej and N.-C. Nguyen, {\it Complex Hessian measures with respect to a background Hermitian form}, https://arxiv.org/abs/2308.10405, to appear in Analysis and PDE.
	
	\bibitem[Li04]{Li04} S.Y. Li, {\it On the Dirichlet problems for symmetric function equations of the eigenvalues of the complex Hessian}, Asian J. Math., {\bf 8} (2004), 87-106.
	
\bibitem[Mi82]{Mi82} M. L. Michelsohn, {\it On the existence of special metrics in complex geometry,} Acta Math. {\bf 149} (1982), no.~3-4, 261--295.
		
\bibitem[Xi96]{Xi96} Y.Xing, {\it Continuity of the complex Monge-Amp\`ere operator}, Proced. AMS, {\bf 124}(1996), no. 2,  457-467.	
		
	
	\end{thebibliography}
\end{document}